\documentclass[12pt]{article}
 \usepackage{amssymb}
 \usepackage{amsmath}
 \usepackage[latin1]{inputenc}
 \usepackage[T1]{fontenc}
\usepackage{times}
\usepackage[colorlinks=true,linkcolor=red,citecolor=blue]{hyperref}
 \usepackage{dsfont}
 \usepackage{amsfonts}
 \usepackage{natbib}
 \usepackage{mathrsfs}
\date{}

 \usepackage{graphicx}

\newtheorem{e-proposition}[theorem]{Proposition}

\newtheorem{e-definition}[theorem]{Definition\rm}

\begin{document}
 \begin{center}
 {\LARGE {On Bayesian Estimation via Divergences}}
 \end{center}
 \begin{center}
 \large{Mohamed Cherfi}\vspace{3mm}

 {\small \it Laboratoire de Statistique Th\'eorique et Appliqu\'ee (LSTA)\\ Equipe d'Accueil 3124\\ Universit\'e Pierre et Marie Curie -- Paris 6\\ Tour 15-25, 2\`eme \'etage\\ 4 place Jussieu\\  75252 Paris cedex 05 }
 \end{center}
\renewcommand{\thefootnote}{}
\footnote{ \vspace*{-4mm} \noindent{{\it E-mail address:}
mohamed.cherfi@gmail.com} }

 \vspace{3mm} \hrule \vspace{3mm} {\small \noindent{\bf Abstract.}
\\
In this Note we introduce a new methodology for Bayesian inference through the use of $\phi$-divergences and the
duality technique. The asymptotic laws of the estimates are established.

\vspace{0.2cm}
 
 \vspace{4mm}\hrule
\section{Introduction}
\noindent Bayesian techniques are particularly attractive since they can incorporate information other than the data into the model in the form of prior distributions. Another feature which make them increasingly attractive is that they can handle models that are difficult to estimate with classical methods by use of simulation techniques, see for instance \cite{Robert2001}.

\noindent The aim of this Note
is to discuss the use of divergences as a basis for Bayesian inference. The use of divergence measures in a Bayesian context has been considered 
in \cite{DeyBirmiwal1994} and \cite{PengDey1995}. Most of this work is concerned with the 
use of divergence measures to study Bayesian robustness when the priors are contaminated and to 
diagnose the effect of outliers.

\noindent In order to estimate model parameters and circumvent possible difficulties encountered with the likelihood function, we follow up common robustification ideas, see for instance \cite{Hanousek1990,Hanousek1994},  and  propose to replace the likelihood in the formula of the posterior distribution by the dual form of the divergence that lead to estimators that are both robust and efficient and include  the expected a posteriori estimator (EAP) as a benchmark. A major advantage of the method is that it does not require
additional accessories such as kernel density estimation or other forms of nonparametric
smoothing to produce nonparametric density estimates of the true underlying density function in contrast with the method proposed by \cite{HookerVidyashankar2011} which is based on the concept of a minimum disparity procedure introduced by \cite{Lindsay1994}.
The plug-in of the empirical distribution function is sufficient for the purpose of estimating the divergence
in the case of i.i.d. data. The proposed estimators are based on integration rather
than optimization. Other reasons, which are commonly put forward to use the proposed approach is computational attractiveness through the use of MCMC and can easily handle a large number of parameters.

\noindent The outline of the Note is as follows. Together with a brief review of definitions and
properties of divergences, Section \ref{estimation} discusses the procedure to obtain the estimates. In Section \ref{Asymp}, we give the limit laws of the proposed estimators. Some final remarks conclude the Note.

\section{Estimation}\label{estimation}
\subsection{Background on dual divergences inference}
\noindent \cite{Keziou2003} and \cite{BroniatowskiKeziou2009} introduced the class of dual divergences estimators for general parametric models,   In the following, we shortly recall their context and definition.

\noindent Recall that the $\phi$-divergence between a bounded
signed measure $Q$ and a probability $P$ on
$\mathscr{D}$, when $Q$ is absolutely continuous with
respect to $P$, is defined by
$$D_\phi(Q,P):=\int_{\mathscr{D}}
\phi\left(\frac{\mathrm{d}Q}{\mathrm{d}P}(x)\right)~\mathrm{d}P(x),$$
where $\phi$ is a convex function from $]-\infty,\infty[$ to
$[0,\infty]$ with $\phi(1)=0$.

\noindent Different choices for $\phi$ have been proposed in the literature. For a good overview, see \cite{Pardo2006}. Well-known class of divergences is the class of the so called
``power divergences'' introduced in \cite{CressieRead1984} (see
also \cite{LieseVajda1987} chapter 2); it contains the most known and used  divergences. They are defined through
the class of convex functions
\begin{equation}  \label{powerdivergence}
x\in ]0,+\infty[ \mapsto\phi_{\gamma}(x):=\frac{x^{\gamma}-\gamma
x+\gamma-1}{\gamma(\gamma-1)}
\end{equation}
if $\gamma\in\mathbb{R}\setminus \left\{0,1\right\}$,
$\phi_{0}(x):=-\log x+x-1$ and $\phi_{1}(x):=x\log x-x+1$.

\noindent Let $X_{1}, \dots, X_{n}$ be an i.i.d. sample
 with p.m. $\mathbb{P}_{\theta_{0}}$. Consider the problem of estimating the population parameters of interest $\theta_{0}$, when the underlying identifiable model is given by $\{\mathbb{P}_{\theta}: \theta\in \Theta\}$ with $\Theta$ a subset of $\mathbb{R}^{d}$. Here the attention is restricted to the case where the probability measures  $\mathbb{P}_{\theta}$ are absolutely
continuous with respect to the same $\sigma$-finite measure $\lambda$; correspondent densities are denoted $p_{\theta}$.

\noindent Let $\phi$ be a function of class $\mathcal{C}^2$, strictly convex and satisfies
\begin{equation}
\int \left| \phi ^{\prime }\left( \frac{p_{\theta }(x)}{p_{\alpha
}(x)}\right) \right|p_{\theta }(x)~\mathrm{d}x<\infty .  \label{condition
integrabilite}
\end{equation}
By Lemma 3.2 in \cite{BroniatowskiKeziou2006}, if the function $\phi$ satisfies:
 There exists $0 <\eta < 1$ such that for all $c$ in $\left[1-\eta, 1 + \eta\right]$,
we can find numbers $c_1$, $c_2$, $c_3$ such that
\begin{equation}\label{EqRem1}
\phi(cx)\leq c_1\phi(x) + c_2 \left|x\right| + c_3, \textrm{ for all real } x,
\end{equation}
then the assumption (\ref{condition
integrabilite}) is satisfied whenever $D_{\phi}(\mathbb{P}_{\theta}, \mathbb{P}_{\alpha})$ is finite. From now on, $\mathcal{U}$ will be the set of $\theta$ and $\alpha$ such that $D_{\phi}(\mathbb{P}_{\theta}, \mathbb{P}_{\alpha})<\infty$. Note that all the real convex functions $\phi_{\gamma}$ pertaining to the class of power divergences defined in (\ref{powerdivergence}) satisfy the condition (\ref{EqRem1}).

\noindent Under (\ref{condition integrabilite}), using Fenchel duality technique, the divergence $D_{\phi}(\theta,\theta_0)$ can be represented as resulting from an optimization procedure, this elegant result was proven in \cite{Keziou2003}, \cite{LieseVajda2006} and \cite{BroniatowskiKeziou2009}. \cite{BroniatowskiKeziou2006} called it the dual form of a divergence, due to its connection with convex analysis.

\noindent Under the above conditions, the $\phi$-divergence:
\begin{equation*}
    D_{\phi}(\mathbb{P}_{\theta}, \mathbb{P}_{\theta_{0}})=\int
\phi\left(\frac{p_{\theta}(x)}{p_{\theta_0}(x)}\right)p_{\theta_0}(x)~\mathrm{d}x,
\end{equation*}
can be represented as the following form:
\begin{equation}\label{Dualrepresentation}
 D_{\phi}(\mathbb{P}_{\theta}, \mathbb{P}_{\theta_{0}})=\sup_{\alpha\in
\mathcal{U}}\int h(\theta,\alpha)~\mathrm{d}\mathbb{P}_{\theta_{0}},
\end{equation}
where $h(\theta,\alpha):x\mapsto h(\theta,\alpha,x)$ and
\begin{equation}\label{Definition-h}
h(\theta,\alpha,x):=\int \phi ^{\prime }\left( \frac{p_{\theta }}{p_{\alpha
}}\right) ~p_{\theta }-\left[ \frac{p_{\theta }(x)}{p_{\alpha
}(x)}\phi ^{\prime }\left( \frac{p_{\theta }(x)}{p_{\alpha }(x)}\right)
-\phi \left( \frac{p_{\theta } (x)}{p_{\alpha }(x)}\right) \right].
\end{equation}
Since the supremum in (\ref{Dualrepresentation}) is unique and is attained in
$\alpha=\theta_0$, independently upon the value of $\theta$, by replacing
the hypothetical probability measure $\mathbb{P}_{\theta_{0}}$
by the empirical measure $\mathbb{P}_{n}$ define the class of estimators of $\theta_{0}$ by
\begin{equation}\label{dualestimator}
\widehat{\alpha}_{\phi}(\theta):=\arg\sup_{\alpha\in \mathcal{U}}\int
h(\theta,\alpha)\mathrm{d}\mathbb{P}_{n},\;\;\theta\in \Theta,
\end{equation}
where $h(\theta,\alpha)$ is the function defined in (\ref{Definition-h}). This class is called ``dual
$\phi$-divergence estimators'' (D$\phi$DE's), see for instance \cite{Keziou2003} and \cite{BroniatowskiKeziou2009}.

\noindent Formula (\ref{dualestimator}) defines a family of $M$-estimators indexed by the function $\phi$ specifying the divergence and by some instrumental value of the parameter $\theta$, called here escort parameter, see also \cite{BroniatowskiVajda2009}.

\noindent Application of dual representation of $\phi$-divergences have been considered by many authors, we cite among others, \cite{KeziouLeoni2008} for semi-parametric two-sample density ratio models, robust tests based on saddlepoint approximations in \cite{TomaLeoni-Aubin2009}, \cite{TomaBroniatowski2010} have proved that this class contains robust and efficient estimators and proposed robust test statistics based on divergences estimators. Bootstrapped $\phi$-divergences estimates are considered in \cite{BouzebdaCherfi2011}, extension of dual $\phi$-divergences estimators to right censored data are introduced in \cite{Cherfi2011a}, for estimation and tests in copula models we refer to \cite{BouzebdaKeziou2010b} and the references therein.
Performances of dual
$\phi$-divergence estimators  for normal models are studied in \cite{Cherfi2011b}.

\subsection{Estimation}
\noindent Let us now turn to the estimation using divergences in our setting.  For the
parameter $\theta$ consider a prior density $\pi(\theta)$ on $\Theta$, and let $\rho(x,\theta) :~\mathbb{R}\times\Theta$ be a
suitable function. Then \cite{Hanousek1990} considered the following {B}ayes-type or {B}-estimator  of $\theta_0$, corresponding to the prior density $\pi(\theta)$ and generated by the function $\rho(x,\theta)$,
\begin{equation}
\widehat{\theta}^*_n=\frac{\int_{\Theta}\theta\exp{\left\lbrace-\sum_{i=1}^n\rho(X_i,\theta)\right\rbrace }\pi(\theta)~\mathrm{d}\theta}{\int_{\Theta}\exp{\left\lbrace-\sum_{i=1}^n\rho(X_i,\theta)\right\rbrace }\pi(\theta)~\mathrm{d}\theta}
\end{equation}
if both integrals exist. This type of estimators is often called Laplace type estimators see for instance \cite{ChernozhukovHong2003}.

\noindent The posterior $M$-estimator is defined as 
\begin{equation}
\widehat{\theta}_n^+=\arg\max_{\theta\in\Theta} \left(-\sum_{i=1}^n\rho(X_i,\theta)+\ln \pi(\theta) \right) .
\end{equation}
\cite{Hanousek1990} showed that $\widehat{\theta}^*_n$ is asymptotically equivalent to the $M$-estimator generated by $\rho$ for a large class of priors and under some conditions on $\rho$ and $\mathbb{P}_{\theta_0}$. The asymptotic equivalence provides the access to the study of
asymptotics for {B}-estimators via the $M$-estimators.

\noindent In the context of the Bayesian methods examined in this Note, instead of a likelihood function, our work will use a criterion function $\displaystyle{\int
h(\theta,\alpha)\mathrm{d}\mathbb{P}_{n}}$.

\noindent Inference is based on the 
$\phi$-posterior 
\begin{equation}
p_{\phi,n}(\alpha|X_1,\cdots,X_n)=\frac{\exp{\left\lbrace n\mathbb{P}_nh(\theta,\alpha)\right\rbrace }\pi(\alpha)}{\int_{\mathcal{U}} \exp{\left\lbrace n\mathbb{P}_nh(\theta,\alpha)\right\rbrace }\pi(\alpha)~\mathrm{d}\alpha}.
\end{equation}
A risk function is the expected loss or error in which the researcher incurs
when choosing a certain value for the parameter estimate.
Let $\mathcal{L}_n(u)$ be a loss function. The risk function takes the form
\begin{equation}
\mathcal{R}_n(\widetilde{\alpha})=\int_{\mathcal{U}}\mathcal{L}_n(\alpha-\widetilde{\alpha})~ p_{\phi,n}(\alpha|X_1,\cdots,X_n)~\mathrm{d}\alpha,
\end{equation}
where $p_{\phi,n}(\alpha|X_1,\cdots,X_n)$ is the $\phi$-posterior density, $\widetilde{\alpha}$ is the selected value, and $\alpha$ is all other possible values we are integrating over. The loss function can penalize the selection of $\alpha$ asymmetrically, and is a function of the selected value and the rest of the possible values of the parameters in $\mathcal{U}$. 

\noindent The dual $\phi$-divergence Bayes type estimator minimizes the expected loss for different forms of the loss
function
\begin{equation}\label{GBdualestimator}
\widehat{\alpha}^*_{\phi}(\theta)=\arg\inf_{\widetilde{\alpha}\in \mathcal{U}}\mathcal{R}_n(\widetilde{\alpha}).
\end{equation}
Choosing different loss functions will change the objective function such
that the estimators bear different interpretations. For instance, when the loss is squared error ($\mathcal{L}_n(u)=|\sqrt{n}u|^2$), for fixed $\theta$, the dual $\phi$-divergence Bayes type estimator is defined as
\begin{equation}\label{Bdualestimator}
\widehat{\alpha}^*_{\phi}(\theta)=\int_{\mathcal{U}} \alpha~ p_{\phi,n}(\alpha|X_1,\cdots,X_n)~\mathrm{d}\alpha:=\frac{\int_{\mathcal{U}} \alpha\exp{\left\lbrace n\mathbb{P}_nh(\theta,\alpha)\right\rbrace }\pi(\alpha)}{\int_{\mathcal{U}} \exp{\left\lbrace n\mathbb{P}_nh(\theta,\alpha)\right\rbrace }\pi(\alpha)~\mathrm{d}\alpha},
\end{equation}
if both integrals exist, other familiar forms obtained for different loss functions are modes, medians
and quantiles. 

\noindent The posterior dual $\phi$-divergences estimator is defined as 
\begin{equation}
\widehat{\alpha}^+_{\phi}(\theta)=\arg\sup_{\alpha\in\mathcal{U}} \left(\mathbb{P}_nh(\theta,\alpha)+\ln \pi(\alpha) \right). 
\end{equation}
It is obvious that posterior dual $\phi$-divergences estimates naturally inherit the properties of dual $\phi$-divergences estimates and hence we focus on dual $\phi$-divergences Bayes type estimators only. 

\newtheorem{Rem1}{Remark}
\begin{Rem1}\label{Rem1}{\rm
\begin{enumerate}
\item The EAP estimator, which is the mean of the posterior distribution, belongs to the
class of estimates (\ref{Bdualestimator}). Indeed, it is obtained when
$\phi (x)=-\log x+x-1$, that is as the dual modified
$ KL_{m}$-divergence estimate. Observe that 
$\displaystyle{\phi^\prime(x)=-\frac{1}{x}+1}$
and 
$x  \phi^\prime(x)-\phi(x)=\log x,$ hence $$\int h(\theta ,\alpha ){\rm{d}}\mathbb{P}_{n}=-\int \log
\left( \frac{{\rm{d}}\mathbb{P}_{\theta }}{{\rm{d}}\mathbb{P}_{\alpha }}\right){\rm{d}}\mathbb{P}_{n}.$$
Keeping in mind definitions (\ref{Bdualestimator}), we get
\begin{equation*}
\widehat{\alpha}^*_{KL_m}(\theta):=\frac{\displaystyle{\int_{\mathcal{U}}\alpha\prod_{i=1}^np_{\alpha}(X_i)\pi(\alpha)~\mathrm{d}\alpha}}{\displaystyle{\int_{\mathcal{U}}\prod_{i=1}^np_{\alpha}(X_i)\pi(\alpha)~\mathrm{d}\alpha}},
\end{equation*}
independently upon $\theta$

\item If new data $X_{n+1},\ldots,X_N$ are obtained, the posterior for the combined data $X_{1},\ldots,X_N$ can be obtained by using posterior after $n$ observations, $p_{\phi,n}(\alpha|X_1,\cdots,X_n)$ as a prior $\alpha$:
\begin{equation*}
p_{\phi,n}(\alpha|X_1,\cdots,X_N)\propto p_{\phi,n}(\alpha|X_1,\cdots,X_n) \times p_{\phi,n}(X_{n+1},\cdots,X_N|\alpha).
\end{equation*}
\end{enumerate}
}
\end{Rem1}
\section{Asymptotic properties}\label{Asymp}
\noindent In this section we state the asymptotic normality of the estimates based on the $\phi$-posterior and evaluate their limiting variance.
The hypotheses handled here are similar to those used in \cite{Keziou2003} and \cite{BroniatowskiKeziou2009} in the frequentist case, these conditions are mild and can be satisfied in most of circumstances. From now on, $\overset{D}{\longrightarrow}$ denotes the convergence in distribution.
\begin{enumerate}
\item[(R.1)] $$\sup_{\alpha\in\Theta}\left\lvert \mathbb{P}_{n}h(\theta,\alpha)-\mathbb{P}_{\theta_{0}}h(\theta,\alpha)\right\rvert\overset{a.s.}{\longrightarrow}0.$$
\item [(R.2)] There exists a neighborhood $N(\theta_0)$
of $\theta_0$ such that the first and second order partial
derivatives (w.r.t $\alpha$) of $\displaystyle{\phi'\left(\frac{p_{\theta} (x)}{p_{\alpha} (x)}\right)p_{\theta} (x)}$
are dominated on $N(\theta_0)$ by some integrable
functions. The third order partial derivatives (w.r.t $\alpha$) of
$h(\theta,\alpha,x)$ are dominated on $N(\theta_0)$ by some
$\mathbb{P}_{\theta_0}$-integrable functions.
\end{enumerate}
Let 
$$S:=-\mathbb{P}_{\theta_0}\frac{\partial^2}{\partial\alpha^2}h\left(\theta,\theta_0\right)\textrm{ and }V:=\mathbb{P}_{\theta_0}\frac{\partial}{\partial\alpha}h\left(\theta,\theta_0\right)^\top\frac{\partial}{\partial\alpha}h\left(\theta,\theta_0\right).$$
Observe that the matrix $S$ is symmetric and positive since the second derivative $\phi''$ is nonnegative by the convexity of $\phi$.
 \begin{enumerate}
 \item[(R.3)]The matrices
$S$ and $V$ are non singular.
\end{enumerate}
For $\alpha$ in an open neighborhood of $\theta_0$, using {\rm{(R.2)}} by a Taylor expansion
\begin{equation}\label{Taylor}
\mathbb{P}_nh\left(\theta,\alpha\right)-\mathbb{P}_nh\left(\theta,\theta_0\right)=(\alpha-\theta_0)^\top U_{n}(\theta_{0})-\dfrac{1}{2}(\alpha-\theta_0)^\top S (\alpha-\theta_0)+R_n(\alpha),  
\end{equation}
 \begin{enumerate}
 \item[(R.4)]Given any $\epsilon> 0$, there exists  $\delta> 0$ such that, the probability
of the event
 \begin{equation}\label{Rn}
 \sup_{|\alpha-\theta_0|\leq\delta}\left|R_n(\alpha)\right|\geq\epsilon
 \end{equation} 
 tends to zero as $n\longrightarrow\infty$.
\end{enumerate}
\newtheorem{Rem2}[Rem1]{Remark}
\begin{Rem2}\label{Rem2}{\rm
\begin{enumerate}
\item Using Example 19.8 in \cite{vanderVaart1998}, it is clear that the class of functions $\displaystyle{\left\{\alpha\mapsto
h(\theta,\alpha);~\alpha\in\Theta\right\}}$ is a
Glivenko-Cantelli class of functions for all fixed $\theta$, that  {\rm{(R.1)}} holds.

\item Conditions {\rm{(R.2)}} and {\rm{(R.3)}} are about usual regularity properties of the underlying model, they guarantee that we can interchange
integration and differentiation and the existence of the variance-covariance matrices, they are similar to regularity conditions used in \cite{Keziou2003} and \cite{BroniatowskiKeziou2009} in the frequentist case.

\item Condition {\rm{(R.4)}} easily holds when there is enough smoothness. It requires that the remainder
term of the expansion can be controlled in a particular way over a neighborhood of $\theta_0$.
\end{enumerate}}
\end{Rem2}
Define
\begin{equation}
t:=\sqrt{n}\left(\alpha-\Delta_n\right),~\Delta_n:=\theta_0+S^{-1}U_{n}(\theta _{0}),
\end{equation}
and $p^{*}_{\phi,n}(t)$ be the $\phi$-posterior density of $t$.

The following theorem states that under some regularity conditions, for large $n$, $p^{*}_{\phi,n}(\cdot)$ is approximately
a random normal density in the $L_1$ sense.
\newtheorem{Theo1}{Theorem}
\begin{Theo1}\label{Theo1}{\rm Let $\pi(\theta)$ be any prior that is continuous and positive at $\theta_{0}$ with $\displaystyle{\int |\theta|\pi(\theta)~\mathrm{d}\theta}$. Then under Conditions  (R.1-4)
\begin{equation}\label{EqTheo1}
\int\left\lvert p^{*}_{\phi,n}(t)-\left(\dfrac{\det S}{2\pi}\right)^{d/2}\exp{\left\lbrace-\frac{1}{2}t^\top S t\right\rbrace }\right\rvert~\mathrm{d}t\overset{P}{\longrightarrow}0.
\end{equation}
}
\end{Theo1}
We now state the principal result of this section. Theorem \ref{Theo2} is concerned with the efficiency and asymptotic
normality of the proposed estimates. See \cite{IbragimovHasminskii1981} and  \cite{Strasser1981} for more on the consistency
and efficiency of Bayes estimators.
\newtheorem{Theo2}[Theo1]{Theorem}
 \begin{Theo2}\label{Theo2}{\rm Let $\pi(\theta)$ be any prior that is continuous and positive at $\theta_{0}$ with $\displaystyle{\int |\theta|\pi(\theta)~\mathrm{d}\theta}$. Assume that Conditions (R.1-4) hold, then as $n$ tends to infinity$$V^{-1/2}S\sqrt{n}\left(\widehat{\alpha}^*_\phi(\theta)-\theta_{0}\right)\overset{d}{\longrightarrow}\mathcal{N}\left(0,I\right).$$
}
 \end{Theo2}
\newtheorem{Rem3}[Rem1]{Remark}
\begin{Rem3}\label{Rem3}{\rm
If $\theta=\theta_{0}$, then $S^\top V^{-1}S=I_{\theta_{0}}$ the information matrix, so that $\widehat{\alpha}^*_\phi(\theta_{0})$ is consistent and asymptotically efficient. The consequence is that the value of the escort parameter should be taken as a consistent estimator of $\theta_{0}$, see \cite{Cherfi2011a,Cherfi2011b} for relevant discussion on this subject.
}
\end{Rem3}
\section{Concluding remarks}
\noindent We have introduced a new estimation procedure in parametric models that combine divergences method with Bayesian analysis, it generalizes the expected a posteriori estimate. The proposed estimators are based on integration rather
than optimization. These estimators are often much easier to compute in practice than
the $\arg\sup$ estimators (\ref{dualestimator}), especially in the high-dimensional setting; see, for example,
the discussion in \cite{LiuTianWei2008}.

\noindent In order to compute these estimators, using Markov chain Monte Carlo methods, we
can draw a Markov chain,
\begin{equation}
\mathbb{S}= (\alpha^{(1)};\alpha^{(2)};\cdots;\alpha^{(B)});
\end{equation}
whose marginal density is approximately given by $p_{\phi,n}(\cdot)$, the $\phi$-posterior distribution.
Then the estimate $\widehat{\alpha}^*_{\phi}(\theta)$ is computed as
\begin{equation}
\widehat{\alpha}^*_{\phi}(\theta)=\frac{1}{B}\sum_{i=1}^B\alpha^{(i)}.
\end{equation}
Consider the construction of confidence intervals for the quantity $f(\theta_0)$, for a given continuously differentiable function $f:\Theta\longrightarrow\mathbb{R}$. Define
\begin{equation}
C_n(\epsilon) :=\inf \left\{x :\int_{f(\alpha)\leq x} \alpha~ p_{\phi,n}(\alpha) ~\mathrm{d}\alpha\geq\epsilon\right\}.
\end{equation}
Then the dual $\phi$-divergence Bayes type estimator confidence interval is given by $\displaystyle{\left[C_n\left(\dfrac{\epsilon}{2}\right); C_n\left(1-\dfrac{\epsilon}{2}\right)\right]}$. These confidence intervals can be constructed simply by taking the $\dfrac{\epsilon}{2}$th and $\dfrac{\epsilon}{2}$th quantiles
of the MCMC sequence
\begin{equation}
f(\mathbb{S})= \left(f(\alpha^{(1)});f(\alpha^{(2)});\cdots;f(\alpha^{(B)})\right),
\end{equation}
and thus are quite simple in practice. 

\noindent The very peculiar choice of the escort parameter defined through $\theta=\theta_0$ has same limit properties as the posterior mean. This result is of some relevance, since it leaves open the choice of the divergence, while keeping good asymptotic properties, we expect that it can also be used directly to provide robust inference, we leave this study for a subsequent paper.

\noindent The problem of the choice of the divergence remain an open question and need more investigation.
\section{Proofs}
Our arguments follow those presented by \cite{LehmannCasella1998}, the main difference is due to  the non-likelihood
setting. See also \cite{ChernozhukovHong2003} for similar arguments.
We often use $M$ to denote a generic finite constant and $I$ to denote the identity matrix.
The smallest-eigenvalue of a matrix $S$ is denoted as ${\rm{mineig}}(S)$.
\subsection{Proof of Theorem \ref{Theo1}}  
Define
\begin{equation}
t:=\sqrt{n}\left(\alpha-\Delta_n\right),~\Delta_n:=\theta_0+S^{-1}U_{n}(\theta _{0}),
\end{equation}
then
\begin{eqnarray}
p^{*}_{\phi,n}(t)&=&\frac{1}{\sqrt{n}}p_{\phi,n}\left(\frac{t}{\sqrt{n}}+\Delta_n\right)\\
\nonumber&=&\dfrac{\pi\left(\dfrac{t}{\sqrt{n}}+\Delta_n\right)\exp{\left\lbrace n\mathbb{P}_nh\left(\theta,\dfrac{t}{\sqrt{n}}+\Delta_n\right)\right\rbrace }}{\displaystyle{\int\pi\left(\dfrac{u}{\sqrt{n}}+\Delta_n\right)\exp{\left\lbrace n\mathbb{P}_nh\left(\theta,\dfrac{u}{\sqrt{n}}+\Delta_n\right)\right\rbrace }~\mathrm{d}u}},\\
&=&\dfrac{\pi\left(\dfrac{t}{\sqrt{n}}+\Delta_n\right)\exp{\left\lbrace\omega(t)\right\rbrace }}{c_n},
\end{eqnarray}
where $$\displaystyle{\omega(t):=n\mathbb{P}_nh\left(\theta,\dfrac{t}{\sqrt{n}}+\Delta_n\right)-n\mathbb{P}_nh\left(\theta,\theta_0\right)-\frac{n}{2}U_{n}(\theta _{0})^\top S^{-1}U_{n}(\theta _{0}),}$$ and $$\displaystyle{c_n:=\int\pi\left(\dfrac{u}{\sqrt{n}}+\Delta_n\right)\exp{\left\lbrace\omega(u)\right\rbrace }~\mathrm{d}u}.$$ 
\newtheorem{Lem1}{Lemma}
\begin{Lem1}\label{Lem1}{\rm
Let 
\begin{equation}\label{J1}
J_1=\int \left\lvert \pi\left(\dfrac{t}{\sqrt{n}}+\Delta_n\right)e^{\omega(t)}-\pi\left(\theta_{0}\right) e^{-\frac{1}{2}t^\top S t}\right\rvert~\mathrm{d}t,
\end{equation}
then if  (R.1-4) hold, $\displaystyle{J_1\overset{P}{\longrightarrow}0}$.}
\end{Lem1} 
By Lemma \ref{Lem1}, we have that
\begin{equation}\label{cn}
c_n\overset{P}{\longrightarrow}\int \pi\left(\theta_{0}\right) e^{-\frac{1}{2}t^\top S t}~\mathrm{d}t=\pi\left(\theta_{0}\right)\sqrt{\dfrac{(2\pi)^{d}}{\lvert\det S\rvert}}.
\end{equation}
Observe that
$$\displaystyle{\int\left\lvert p^{*}_{\phi,n}(t)-\left(\dfrac{\det S}{2\pi}\right)^{d/2}\exp{\left\lbrace-\frac{1}{2}t^\top S t\right\rbrace }\right\rvert~\mathrm{d}t=\dfrac{J}{c_n},}$$
 where
\begin{equation}
J:=\int \left\lvert \pi\left(\dfrac{t}{\sqrt{n}}+\Delta_n\right)e^{\omega(t)}-c_n\sqrt{\dfrac{\lvert\det S\rvert}{(2\pi)^{d}}} e^{-\frac{1}{2}t^\top S t}\right\rvert~\mathrm{d}t
\end{equation}
By (\ref{cn}), to show (\ref{EqTheo1}) it is enough to show that $\displaystyle{J\overset{P}{\longrightarrow}0}$. But, $\displaystyle{J\leq J_1+J_2}$ where $J_1$ is given by (\ref{J1}) and
\begin{equation*}
J_2=\int \left\lvert c_n\sqrt{\dfrac{\lvert\det S\rvert}{(2\pi)^{d}}} e^{-\frac{1}{2}t^\top S t}-\pi\left(\theta_{0}\right) e^{-\frac{1}{2}t^\top S t}\right\rvert~\mathrm{d}t.
\end{equation*}
Observe that 
\begin{equation}\label{J2}
J_2=\left\lvert c_n\sqrt{\dfrac{\lvert\det S\rvert}{(2\pi)^{d}}}-\pi\left(\theta_{0}\right)\right\rvert\int e^{-\frac{1}{2}t^\top S t}~\mathrm{d}t\overset{P}{\longrightarrow}0.
\end{equation}
By Lemma \ref{Lem1} and  (\ref{J2}), $J_1$ and $J_2$ tend to zero in probability, and this completes the proof.
\subsection*{Proof of Lemma \ref{Lem1}}
\noindent Let 
$$U_{n}(\theta _{0}):=P_{n}\frac{\partial }{\partial \alpha}h(\theta ,\theta_0).$$
Using {\rm{(R.2)}} and {\rm{(R.3)}} in connection with the central limit theorem (CLT), we can see that
\begin{equation}\label{CLTUn}
\sqrt{n}V^{-1/2}U_{n}(\theta_{0})\overset{d}{\longrightarrow}\mathcal{N}(0,I).
\end{equation}
Write
\begin{equation}\label{wt}
\omega(t)=-\frac{1}{2}t^\top S t+R_n\left(\dfrac{t}{\sqrt{n}}+\Delta_n\right).
\end{equation}
To prove that the integral (\ref{J1}) tends to zero in probability, divide
the range of integration into the three parts:
\begin{itemize}
\item[(i)] $|t|\leq M$,
\item[(ii)] $|t|\geq \delta \sqrt{n} $,
\item[(iii)] $M<|t|<\delta \sqrt{n}$, 
\end{itemize}
and show that the integral over each of the three tends to zero in
probability.

\noindent Part(i):
\begin{equation*}
\int_{|t|\leq M} \left\lvert \pi\left(\dfrac{t}{\sqrt{n}}+\Delta_n\right)e^{\omega(t)}-\pi\left(\theta_{0}\right) e^{-\frac{1}{2}t^\top S t}\right\rvert~\mathrm{d}t\overset{P}{\longrightarrow}0.
\end{equation*}
To prove this result, we shall show that for every $0 < M <\infty$,
\begin{equation}\label{i}
\sup_{|t|\leq M}\left\lvert \pi\left(\dfrac{t}{\sqrt{n}}+\Delta_n\right)e^{\omega(t)}-\pi\left(\theta_{0}\right) e^{-\frac{1}{2}t^\top S t}\right\rvert\overset{P}{\longrightarrow}0.
\end{equation}
Substituting the
expression (\ref{wt}) for $\omega(t)$, (\ref{i}) is seen to follow from the following two facts
\begin{equation}\label{Fact1}
\sup_{|t|\leq M}\left\lvert \pi\left(\dfrac{t}{\sqrt{n}}+\Delta_n\right)-\pi\left(\theta_{0}\right) \right\rvert\overset{P}{\longrightarrow}0
\end{equation}
and 
\begin{equation}\label{Fact2}
\sup_{|t|\leq M}\left\lvert R_n\left(\dfrac{t}{\sqrt{n}}+\Delta_n\right)\right\rvert\overset{P}{\longrightarrow}0.
\end{equation}
The first fact is obvious from the continuity of $\pi$ and because by Condition {\rm{(R.3)}}  and (\ref{CLTUn}):
\begin{equation}\label{A11}
\sqrt{n}S^{-1}U_{n}(\theta _{0})=O_P(1),
\end{equation}
so that $\displaystyle{\Delta_n \overset{P}{\longrightarrow}\theta_0}$.

\noindent Given (\ref{A11}), the second fact follows from Condition {\rm{(R.2)}}, and 
\begin{equation*}
\sup_{|t|\leq M}\left\lvert \Delta_n+\dfrac{t}{\sqrt{n}}-\theta_0\right\rvert=O_P(\dfrac{1}{\sqrt{n}}).
\end{equation*}
\noindent Part(ii):
\begin{equation*}
\int_{M<|t|<\delta\sqrt{n}} \left\lvert \pi\left(\dfrac{t}{\sqrt{n}}+\Delta_n\right)e^{\omega(t)}-\pi\left(\theta_{0}\right) e^{-\frac{1}{2}t^\top S t}\right\rvert~\mathrm{d}t\overset{P}{\longrightarrow}0.
\end{equation*}
For the second part, since the integral of the second term is finite and can be made arbitrarily small by
setting $M$ large, it suffices to show that for the integrand of the first term is bounded by an integrable function with probability $\geq 1-\epsilon$. More precisely, we shall
show that given $\epsilon> 0$, there exists $ \delta> 0$ and $C <\infty$ such that for sufficiently
large $n$,
\begin{equation}\label{Bound}
P\left[\pi\left(\dfrac{t}{\sqrt{n}}+\Delta_n\right)e^{\omega(t)}\leq Ce^{-\frac{1}{4}t^\top S t} \textrm{ for all }M<|t|<\delta\sqrt{n}\right]\geq 1-\epsilon .
\end{equation}
By the fact that $\displaystyle{\Delta_n \overset{P}{\longrightarrow}\theta_0}$ and the continuity of $\pi$, we can drop the factor $\displaystyle{\pi\left(\dfrac{t}{\sqrt{n}}+\Delta_n\right)}$ from consideration, so that it remains to establish such a bound for $\displaystyle{\exp{\left\{\omega(t)\right\}}}$. By definition of $\omega(t)$ (\ref{wt})
\begin{equation}\label{A12}
\exp{\left\{\omega(t)\right\}}\leq\exp{\left\{-\frac{1}{2}t^\top S t+R_n\left(\dfrac{t}{\sqrt{n}}+\Delta_n\right)\right\}}.
\end{equation}
Since $\displaystyle{\left|\Delta_n-\theta_0\right|=o_P(1)}$, it follows that with probability arbitrarily close to $1$, for $n$ sufficiently
large,
$$\displaystyle{\left|\Delta_n+\dfrac{t}{\sqrt{n}}-\theta_0\right|<2\delta'\textrm{ for all }|t|\leq\delta'\sqrt{n}.}$$
Thus, by Condition {\rm{(R.4)}}, there exists some small $\delta'$ and large $M$ such that the latter inequality implies
\begin{equation*}
P\left[\sup_{M\leq|t|\leq\delta'\sqrt{n}}\left|R_n\left(\dfrac{t}{\sqrt{n}}+\Delta_n\right)\right|\leq \dfrac{1}{4}{\rm{mineig}}(S)\right]\geq 1-\epsilon .
\end{equation*}
Combining this fact with (\ref{A11}), we see that (\ref{A12}), for some $C>0$, is
\begin{equation}
\exp{\left\{\omega(t)\right\}}\leq C\exp{\left\{-\frac{1}{2}t^\top S t\right\}},
\end{equation}
for all $t$ satisfying {\rm(ii)}, with probability arbitrarily close to $1$, and this establishes (\ref{Bound}).

\noindent Part(iii):
\begin{equation*}
\int_{|t|\geq\delta\sqrt{n}} \left\lvert \pi\left(\dfrac{t}{\sqrt{n}}+\Delta_n\right)e^{\omega(t)}-\pi\left(\theta_{0}\right) e^{-\frac{1}{2}t^\top S t}\right\rvert~\mathrm{d}t\overset{P}{\longrightarrow}0.
\end{equation*}
As in (ii), the second term in the integrand can be
neglected. Therefore we only need
to show
\begin{equation*}
\int_{|t|\geq\delta\sqrt{n}} \pi\left(\dfrac{t}{\sqrt{n}}+\Delta_n\right)e^{\omega(t)}~\mathrm{d}t\overset{P}{\longrightarrow}0.
\end{equation*}
Recalling the definition of $t$, the term is bounded by
\begin{equation*}
\int_{|\alpha-\Delta_n|\geq\delta} \pi\left(\alpha\right)\exp{\left\lbrace n\mathbb{P}_nh\left(\theta,\alpha\right)-n\mathbb{P}_nh\left(\theta,\theta_0\right)-\frac{n}{2}U_{n}(\theta _{0})^\top S^{-1}U_{n}(\theta _{0})\right\rbrace }~\mathrm{d}\alpha.
\end{equation*}
By {\rm{(R.4)}}, for any $\delta>0$, there exists $\epsilon>0$, such that
\begin{equation}
\sup_{|\alpha-\theta_0|\geq \delta}\left( \mathbb{P}_nh\left(\theta,\alpha\right)-\mathbb{P}_nh\left(\theta,\theta_0\right)\right)\leq -\epsilon. 
\end{equation}
Since $\displaystyle{\Delta_n \overset{P}{\longrightarrow}\theta_0}$, therefore with probability tending to $1$, there exists $\epsilon$ such that
\begin{equation*}
\sup_{|\alpha-\Delta_n|\geq\delta}\exp{\left\lbrace n\mathbb{P}_nh\left(\theta,\alpha\right)-n\mathbb{P}_nh\left(\theta,\theta_0\right)\right\rbrace }\leq e^{-n\epsilon}.
\end{equation*} 
Since $\displaystyle{\exp{\left\lbrace-\frac{n}{2}U_{n}(\theta _{0})^\top S^{-1}U_{n}(\theta _{0})\right\rbrace }=O_P(1)}$, the entire term is bounded by
\begin{equation*}
C\sqrt{n}e^{-n\epsilon}\int \pi\left(\alpha\right)~\mathrm{d}\alpha=o_P(1),
\end{equation*}
with probability tending to $1$. 

\noindent The entire proof is now completed by combining all terms.
\subsection{Proof of Theorem \ref{Theo2}}  
\noindent We have
$$V^{-1/2}S\sqrt{n}\left(\widehat{\alpha}^*_\phi(\theta)-\theta_{0}\right)=V^{-1/2}S\sqrt{n}\left(\widehat{\alpha}^*_\phi(\theta)-\Delta_n\right)+V^{-1/2}S\sqrt{n}\left(\Delta_n-\theta_{0}\right).$$ 
By the CLT, the second term has the limit distribution $\mathcal{N}\left(0,I\right)$, so that it only remains to show that
\begin{equation}
\sqrt{n}\left(\widehat{\alpha}^*_\phi(\theta)-\Delta_n\right)\overset{P}{\longrightarrow}0.
\end{equation}
Observe that
\begin{eqnarray*}
\widehat{\alpha}^*_{\phi}(\theta)&=&\int\alpha~ p_{\phi,n}(\alpha)~\mathrm{d}\alpha\\
&=&\int\left(\frac{t}{\sqrt{n}}+\Delta_n\right)~ p^{*}_{\phi,n}(t)~\mathrm{d}t\\
&=&\frac{1}{\sqrt{n}}\int t~ p^{*}_{\phi,n}(t)~\mathrm{d}t+\Delta_n,
\end{eqnarray*}
and hence
\begin{equation}
\sqrt{n}\left(\widehat{\alpha}^*_\phi(\theta)-\Delta_n\right)=\int t~ p^{*}_{\phi,n}(t)~\mathrm{d}t.
\end{equation}
Thus,
\begin{eqnarray*}
\sqrt{n}\left|\widehat{\alpha}^*_\phi(\theta)-\Delta_n\right|&=&\left|\int t~ p^{*}_{\phi,n}(t)~\mathrm{d}t-\int t~ \left(\dfrac{\det S}{2\pi}\right)^{d/2}\exp{\left\lbrace-\frac{1}{2}t^\top S t\right\rbrace }~\mathrm{d}t\right|\\
&\leq&\int |t|~ \left|p^{*}_{\phi,n}(t)~\mathrm{d}t-\left(\dfrac{\det S}{2\pi}\right)^{d/2}\exp{\left\lbrace-\frac{1}{2}t^\top S t\right\rbrace }\right|~\mathrm{d}t,
\end{eqnarray*}
which tends to zero in probability by Theorem \ref{Theo1}.
\def\ocirc#1{\ifmmode\setbox0=\hbox{$#1$}\dimen0=\ht0 \advance\dimen0
  by1pt\rlap{\hbox to\wd0{\hss\raise\dimen0
  \hbox{\hskip.2em$\scriptscriptstyle\circ$}\hss}}#1\else {\accent"17 #1}\fi}

\end{document}